\documentclass[11pt,english,reqno]{amsart}

\usepackage{amsthm,amsmath,amssymb,amsfonts}
\usepackage{mathtools,mleftright}
\usepackage{microtype,setspace,epigraph}
\usepackage[dvipsnames]{xcolor}
\usepackage{longtable,booktabs,float,caption,subcaption} 
\usepackage[a4paper,margin=1in]{geometry}

\definecolor{webgreen}{rgb}{0,.5,0}
\definecolor{webbrown}{rgb}{.6,0,0}

\usepackage[
  pdfauthor={},
  pdfkeywords={},
  pdftitle={},
  pdfcreator={},
  pdfproducer={},
  colorlinks,linkcolor=blue,filecolor=webbrown,citecolor=webgreen]{hyperref}

\mleftright

\theoremstyle{plain}
\newtheorem{theorem}{Theorem}[section]

\newtheorem{corollary}[theorem]{Corollary}

\theoremstyle{definition}

\newtheorem{remark}[theorem]{Remark}

\numberwithin{equation}{section}
\numberwithin{table}{section}
\numberwithin{figure}{section}

\newcommand{\ZZ}{\mathbb{Z}}
\newcommand{\QQ}{\mathbb{Q}}
\newcommand{\RR}{\mathbb{R}}
\newcommand{\CC}{\mathbb{C}}

\DeclareMathOperator{\Int}{Int}
\DeclareMathOperator{\Mat}{M}

\newcommand{\AP}{\mathbf{A}}
\newcommand{\AB}{\mathbb{A}}
\newcommand{\BE}{\mathcal{B}}
\newcommand{\BN}{\mathbf{B}}
\newcommand{\BB}{\mathbb{B}}
\newcommand{\EN}{\mathbf{E}}

\newcommand{\HN}{\mathbf{H}}
\newcommand{\IM}{\mathcal{I}}
\newcommand{\FM}{\mathcal{F}}
\newcommand{\FN}{\mathbf{F}}
\newcommand{\LN}{\mathbf{L}}

\newcommand{\ST}{\mathbf{S}}
\newcommand{\RT}{\mathfrak{R}}
\newcommand{\CT}{\mathfrak{C}}

\newcommand{\textq}[1]{\quad \text{#1} \quad}
\newcommand{\andq}{\textq{and}}
\newcommand{\set}[1]{\left\{#1\right\}}
\newcommand{\seq}[1]{\left(#1_\nu\right)_{\nu \geq 0}}
\newcommand{\seqs}[1]{\left(#1\right)_{\nu \geq 0}}
\newcommand{\valueat}[1]{{\,}_{\big|\, #1}}
\newcommand{\phrase}[1]{``#1''}
\newcommand{\ditto}{\texttt{"}}
\newcommand{\rhs}{\mathrm{RHS}}
\newcommand{\seqnum}[1]{\href{https://oeis.org/#1}{\rm \underline{#1}}}

\begin{document}

\title{Boomerang identities and binomial transform}
\author{Bernd C. Kellner}
\address{G\"ottingen, Germany}
\email{bk@bernoulli.org}

\begin{abstract}
In Pascal's triangle, the boomerang identities, which were recently introduced,
combine and extend the well-known hockey-stick identities of the binomial 
coefficients by allowing additional terms. We show that the binomial transform
of sequences in a ring play a central role in this context. Subsequently,
we show several applications to Bell and Bernoulli numbers and polynomials,
Appell polynomials, as well as, to harmonic, Fibonacci, Lucas, and Stirling numbers.
\end{abstract}

\maketitle

\vspace*{-6ex}
\section{Introduction}

The \emph{boomerang} identities in \emph{Pascal's triangle} were recently 
introduced by the author~\cite{Kellner:2026}, which combine and extend the 
well-known \emph{hockey-stick} identities of the binomial coefficients by allowing 
additional terms. Figures~\ref{fig:boom} and~\ref{fig:hock} below show examples 
of shapes of the boomerang and hockey-stick identities, respectively.
In the present paper, we show further properties in connection with sequences 
in rings and their binomial transforms. 

Let $\seq{s}$ be a sequence in a ring $R$. The \emph{binomial transform} of 
$\seq{s}$ is defined by
\vspace*{-0.8ex}
\begin{align} 
  t_n = \sum_{\nu=0}^{n} \binom{n}{\nu} s_\nu, 
  \label{eq:btr1}
\shortintertext{providing the sequence $\seq{t}$. 
Its unique inversion is given by iterated forward differences via}
  s_n = \sum_{\nu=0}^{n} \binom{n}{\nu} (-1)^{n-\nu} t_\nu. 
  \label{eq:btr2}
\end{align}
Usually, these transforms are defined over fields like $\RR$ or $\CC$, 
and they can be differently expressed depending on the choice of signs, 
see Riordan~\cite[Chap.~2,\,pp.~43--45]{Riordan:1968}.
We use the shorter operator notation for \eqref{eq:btr1} and \eqref{eq:btr2}, namely,
\[
  t_n = \nabla^n s \andq s_n = \Delta^{\!n} \, t,
\]
where $\nabla = E+I$ and $\Delta = E-I$ with $E$ and $I$ being the shift and identity 
operator acting on sequences, respectively. The expansions of $\nabla^n$ and $\Delta^{\!n}$ 
then act directly on the underlying abelian group. As a consequence, the ring $R$ does 
not necessarily have to be commutative or have a unit element.
For example, take $R = \Mat_n(\RR)$ $(n \geq 2)$ as the ring of $n \times n$ matrices 
over $\RR$ or $R = 2\ZZ$. Theorem~\ref{thm:main} gives the main result of the paper.

\begin{theorem} \label{thm:main}
Let $\seq{s}$ and $\seq{t}$ be any sequences in a ring $R$ (not necessarily unital).
For $n \geq k \geq 0$, we have the boomerang identities
\begin{align}
  \sum_{\nu=0}^{k} \binom{n+1}{k-\nu} s_\nu
  &= \sum_{\nu=0}^{k} \binom{n-\nu}{k-\nu} t_\nu 
  \label{eq:seq1} 
\shortintertext{and}
  \sum_{\nu=k}^{n} \binom{n+1}{\nu+1} s_{\nu-k} 
  &= \sum_{\nu=k}^{n} \binom{\nu}{k} t_{n-\nu},
  \label{eq:seq2} 
\shortintertext{as well as the limiting case}
  \sum_{\nu=0}^{n} \binom{n+1}{\nu+1} s_\nu
  &= \sum_{\nu=0}^{n} t_\nu, 
  \label{eq:seq3} 
\end{align}
if and only if
\begin{equation} \label{eq:seq4} 
  \seq{t} = \seqs{\nabla^\nu s} \textq{or} \seq{s} = \seqs{\Delta^{\!\nu} t}.
\end{equation}
\end{theorem}

\begin{figure}[H]
\begin{center}
\begin{subfigure}[b]{5cm}
\includegraphics[width=5cm]{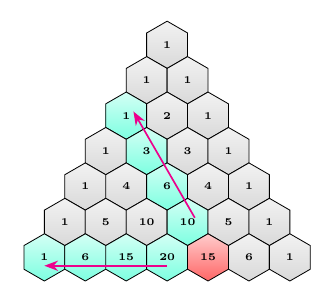}
\vspace*{-5ex}
\caption{}
\end{subfigure}
\begin{subfigure}[b]{5cm}
\includegraphics[width=5cm]{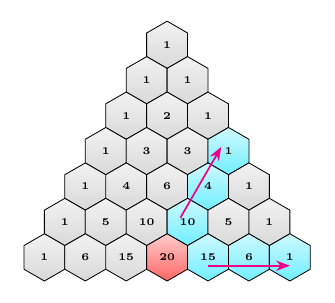}
\vspace*{-5ex}
\caption{}
\end{subfigure}
\begin{subfigure}[b]{5cm}
\includegraphics[width=5cm]{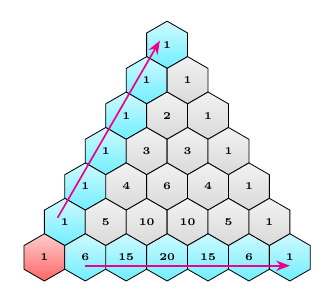}
\vspace*{-5ex}
\caption{}
\end{subfigure}
\caption{Colored shapes of boomerangs in Pascal's triangle, where the figures 
(a)~--~(c) correspond to \eqref{eq:seq1}~--~\eqref{eq:seq3}, respectively.
The red hexagonal areas mark the omitted vertices, from which the arrows point 
in the direction of the increasing indices of the sequences.} 
\label{fig:boom}
\end{center}
\end{figure}

As special cases, the sequences $\seq{s} = (1,0,0,\ldots)$ and $\seq{t} = (1,0,0,\ldots)$
generate the two different kinds of the hockey-stick identities as shown in 
Figure~\ref{fig:hock} below, respectively.

\begin{corollary} \label{cor:hock}
Let $n \geq k \geq 0$. We have the following statements.
(1) The sequence $\seq{s} = (1,0,0,\ldots)$ implies $\seq{t} = (1,1,1,\ldots)$,
giving the hockey-stick identities
\begin{align}
  \sum_{\nu=0}^{k} \binom{n-k+\nu}{\nu} &= \binom{n+1}{k}
  \label{eq:hs1}
\shortintertext{and}
  \sum_{\nu=k}^{n} \binom{\nu}{k} &= \binom{n+1}{k+1}.
  \label{eq:hs2}
\end{align}
(2) The sequence $\seq{s} = (1,-1,1,-1,\ldots)$ implies $\seq{t} = (1,0,0,\ldots)$,
giving the (lying) hockey-stick identities
\begin{align}
  \sum_{\nu=0}^{k} \binom{n+1}{\nu} (-1)^{k-\nu} &= \binom{n}{k} 
  \label{eq:hs3} 
\shortintertext{and}
  \sum_{\nu=k}^{n} \binom{n+1}{\nu+1} (-1)^{\nu-k} &= \binom{n}{k}. 
  \label{eq:hs4}
\end{align}
\end{corollary}

\begin{figure}[H]
\begin{center}
\begin{subfigure}[b]{5cm}
\includegraphics[width=5cm]{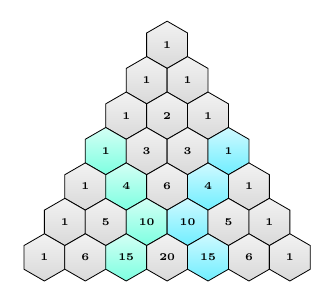}
\vspace*{-5ex}
\caption{}
\end{subfigure}
\begin{subfigure}[b]{5cm}
\includegraphics[width=5cm]{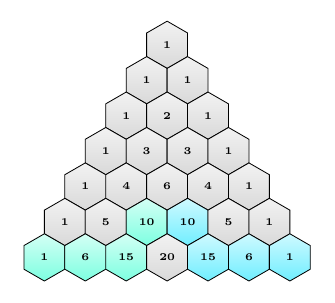}
\vspace*{-5ex}
\caption{}
\end{subfigure}
\caption{Symmetric colored shapes of hockey sticks in Pascal's triangle. 
(a) The shapes on the left and right correspond to 
\eqref{eq:hs1} and \eqref{eq:hs2}, respectively.
(b) The shapes on the left and right correspond to 
\eqref{eq:hs3} and \eqref{eq:hs4}, respectively.} 
\label{fig:hock}
\end{center}
\end{figure}

The results can be applied to functions $f: \RR \to \CC$, where we define
\[
  \Delta^{\!\nu} f = \Delta^{\!\nu} f(x) \valueat{x=0},
\]
indicating the variable and an initial value. The Gregory-Newton expansion of 
a polynomial $p \in \CC[x]$ with $d = \deg p$ is given by
\[
  p(x) = \sum_{\nu=0}^{d} \binom{x}{\nu} \Delta^{\!\nu} p,
\]
where the coefficients $\Delta^{\!\nu} p$ vanish for $\nu > d$. The zero polynomial 
has degree ${d = -\infty}$. A~polynomial $p \in \QQ[x]$ is said to be an 
\emph{integer-valued polynomial} if $p(n) \in \ZZ$ for all $n \in \ZZ$. 
The set of all such polynomials forms a ring denoted by $\Int(\ZZ)$, which is 
a free $\ZZ$-module and the binomial polynomials $\set{\binom{x}{\nu}}_{\nu \geq 0}$ 
provide a basis. 

\begin{theorem} \label{thm:func}
Let $f: \RR \to \CC$. For $n \geq k \geq 0$, we have
\begin{align}
  \sum_{\nu=0}^{k} \binom{n+1}{k-\nu} \Delta^{\!\nu} f
  &= \sum_{\nu=0}^{k} \binom{n-\nu}{k-\nu} f(\nu) 
  \label{eq:fp1} 
\shortintertext{and}
  \sum_{\nu=k}^{n} \binom{n+1}{\nu+1} \Delta^{\!\nu-k} f 
  &= \sum_{\nu=k}^{n} \binom{\nu}{k} f(n-\nu), 
  \label{eq:fp2} 
\shortintertext{as well as the limiting case}
  \sum_{\nu=0}^{n} \binom{n+1}{\nu+1} \Delta^{\!\nu} f
  &= \sum_{\nu=0}^{n} f(\nu). 
  \label{eq:fp3} 
\end{align}
We obtain truncated identities in the following cases.
If $f \in \CC[x]$ is a polynomial with $d = \deg f$, then $\Delta^{\!\nu} f = 0$ 
for $\nu > d$. Moveover, if $f \in \Int(\ZZ)$, then all equations lie in $\ZZ$.
\end{theorem}

\begin{figure}[H]
\begin{center}
\begin{subfigure}[b]{5cm}
\includegraphics[width=5cm]{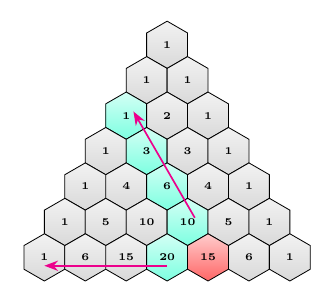}
\vspace*{-5ex}
\caption{}
\end{subfigure}
\begin{subfigure}[b]{5cm}
\includegraphics[width=5cm]{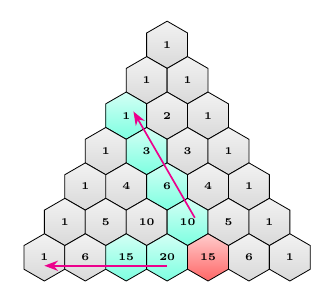}
\vspace*{-5ex}
\caption{}
\end{subfigure}
\begin{subfigure}[b]{5cm}
\includegraphics[width=5cm]{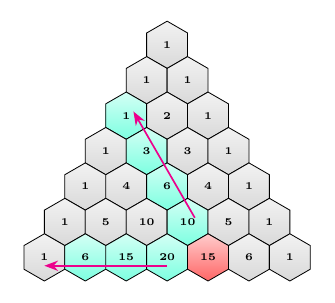}
\vspace*{-5ex}
\caption{}
\end{subfigure}
\caption{Colored shapes of truncated boomerangs in Pascal's triangle. The figures 
(a)~--~(c) correspond to polynomials with degrees 0, 1, 2, respectively.} 
\label{fig:trunc}
\end{center}
\end{figure}

As a special case, we consider binomial polynomials $\binom{x}{\nu} \in \Int(\ZZ)$,
where we usually set ${\binom{x}{\nu} = 0}$ for ${\nu < 0}$.
See Figure~\ref{fig:trunc} for truncated shapes of boomerangs for degrees 
$d = 0, 1, 2$.

\begin{corollary} \label{cor:bin}
Let $d \geq 0$. For $n \geq k \geq 0$, we have
\begin{align*}
  \sum_{\nu=0}^{k} \binom{n+1}{k-\nu} \binom{x}{d-\nu} 
  &= \sum_{\nu=0}^{k} \binom{n-\nu}{k-\nu} \binom{x+\nu}{d} 
\shortintertext{and}
  \sum_{\nu=k}^{n} \binom{n+1}{\nu+1} \binom{x}{d+k-\nu}  
  &= \sum_{\nu=k}^{n} \binom{\nu}{k} \binom{x+n-\nu}{d}, 
\shortintertext{as well as the limiting case}
  \sum_{\nu=0}^{n} \binom{n+1}{\nu+1} \binom{x}{d-\nu} 
  &= \sum_{\nu=0}^{n} \binom{x+\nu}{d}. 
\end{align*}
\end{corollary}

\vspace*{-2.1ex}

\section{Applications}

In this section, we show several applications of the binomial transform in context
of the boomerang identities. The derived identities below are listed in the appendix, 
see Table~\ref{tab:sum}. We omit the corresponding identities that would arise from 
\eqref{eq:seq2} due to symmetry. As usual, empty sums are defined as $0$, and $0^0 = 1$.

\subsection{Bell numbers and polynomials}

The sequence of the Bell numbers for $n \geq 0$ begins
\[
  \BE_n = 1, 1, 2, 5, 15, 52, 203, 877, 4140, 21147, \ldots,
\]
see sequence \seqnum{A000110} in the OEIS~\cite{OEIS}. These numbers satisfy
the recurrence relation
\begin{equation} \label{eq:bn}
  \BE_{n+1} = \sum_{\nu=0}^{n} \binom{n}{\nu} \BE_\nu.
\end{equation}
The Bell polynomials $\BE_n(x)$ (see Bell~\cite{Bell:1934}), also called exponential 
polynomials, as well as, Touchard polynomials (see Touchard~\cite{Touchard:1933}). 
They generalize the Bell numbers by 
\[
  \BE_n(x) = \sum_{\nu=0}^{n} \ST_2(n,\nu) \, x^\nu
\]
with the Stirling numbers of the second kind $\ST_2$. 
Since $\BE_n = \BE_n(1)$, we have 
\[
  \BE_n = \sum_{\nu=0}^{n} \ST_2(n,\nu).
\]
The recurrence relation reads
\begin{equation} \label{eq:bnp}
  \BE_{n+1}(x) = x \sum_{\nu=0}^{n} \binom{n}{\nu} \BE_\nu(x).
\end{equation}

Since the recurrences \eqref{eq:bn} and \eqref{eq:bnp} describe binomial transforms, 
we deduce from Theorem~\ref{thm:main} with $R=\ZZ[x]$ the following corollary.

\begin{corollary} 
Let $n \geq k \geq 0$, then
\begin{align*}
  \sum_{\nu=0}^{k} \binom{n+1}{k-\nu} x \, \BE_\nu(x)
  &= \sum_{\nu=0}^{k} \binom{n-\nu}{k-\nu} \BE_{\nu+1}(x)
\shortintertext{and}
  \sum_{\nu=0}^{k} \binom{n+1}{k-\nu} \BE_\nu
  &= \sum_{\nu=0}^{k} \binom{n-\nu}{k-\nu} \BE_{\nu+1}.
\shortintertext{In particular,}
  \sum_{\nu=0}^{n} \binom{n+1}{\nu+1} x \, \BE_\nu(x)
  &= \sum_{\nu=0}^{n} \BE_{\nu+1}(x)
\shortintertext{and}
  \sum_{\nu=0}^{n} \binom{n+1}{\nu+1} \BE_\nu
  &= \sum_{\nu=0}^{n} \BE_{\nu+1}. 
\end{align*}
\end{corollary}

\subsection{Stirling numbers of the second kind}

The numbers $\ST_2(n,k)$ satisfy
\begin{equation} \label{eq:stdef}
  \ST_2(n,k) = \frac{1}{k!} \Delta^{\!k} \, \nu^n \valueat{\nu=0} 
\end{equation}
and vanish for $k > n$. By Theorem~\ref{thm:func}, we obtain the following result.

\begin{corollary} 
Let $m \geq 0$. For $n \geq k \geq 0$, we have
\begin{align}
  \sum_{\nu=0}^{k} \binom{n+1}{k-\nu} \nu! \, \ST_2(m,\nu)
  &= \sum_{\nu=0}^{k} \binom{n-\nu}{k-\nu} \nu^m.
  \nonumber
\shortintertext{In particular,}
  \sum_{\nu=0}^{n} \binom{n+1}{\nu+1} \nu! \, \ST_2(m,\nu)
  &= \sum_{\nu=0}^{n} \nu^m. 
  \label{eq:psum}
\end{align}
\end{corollary}

Formula \eqref{eq:psum} implies the well-known representation of a power sum 
\[
  S_m(n) = \sum_{\nu=0}^{n-1} \nu^m
\]
with the Stirling numbers as a polynomial in $\Int(\ZZ)$. Due to \eqref{eq:stdef},
we get 
\[
  S_m(x) = \sum_{\nu=0}^{m} \binom{x}{\nu+1} \nu! \, \ST_2(m,\nu)
\]
with $\deg S_m = m+1$.

\subsection{Harmonic numbers}

Define the harmonic numbers for $n \geq 1$ by 
\[
  \HN_n = \sum_{\nu=1}^{n} \frac{1}{\nu}. 
\]
We have the binomial transform (cf.\ Graham et al.~\cite[Sec.\,6.4, pp.\,281--282]{GKP:1994})
\[
  \HN_n = \sum_{\nu=1}^{n} \binom{n}{\nu} (-1)^{\nu+1} \frac{1}{\nu},
\]
implying the following result.

\begin{corollary} 
Let $n \geq k \geq 0$, then
\begin{align*}
  \sum_{\nu=1}^{k} \binom{n+1}{k-\nu} (-1)^{\nu+1} \frac{1}{\nu}
  &= \sum_{\nu=1}^{k} \binom{n-\nu}{k-\nu} \HN_\nu.
\shortintertext{In particular,}
  \sum_{\nu=1}^{n} \binom{n+1}{\nu+1} (-1)^{\nu+1} \frac{1}{\nu}
  &= \sum_{\nu=1}^{n} \HN_\nu. 
\end{align*}
\end{corollary}

\subsection{Bernoulli numbers and polynomials}

The Bernoulli polynomials are defined by the exponential generating function
(cf.\ Comtet~\cite[Sec.\,1.14, pp.\,48--49]{Comtet:1974}) 
\[
  \frac{t e^{xt}}{e^t - 1} = \sum_{n \geq 0} \BN_n(x) \frac{t^n}{n!} 
  \quad (|t| < 2\pi), 
\]
where $\BN_n = \BN_n(0) \in \QQ$ are the Bernoulli numbers. 
For $r,\ell \geq 0$, define the numbers
\[
  \BB_{r,\ell} = \sum_{\nu=0}^{r} \binom{r}{\nu} \BN_{\ell+\nu},
\]
arising from the binomial transform. They obey the well-known reciprocity
relation
\begin{equation} \label{eq:bnid}
  (-1)^r \, \BB_{r,\ell} = (-1)^\ell \, \BB_{\ell,r}.
\end{equation}
The above identity has been recovered by several authors; 
see Kellner~\cite[Sec.\,3]{Kellner:2022} for historical notes, 
see also Prodinger~\cite{Prodinger:2014} for a short proof of \eqref{eq:bnid}.
For example, Carlitz~\cite{Carlitz:1971} gave \eqref{eq:bnid} as a problem in 1971. 
However, we note as in \cite{Kellner:2022} that identity~\eqref{eq:bnid} goes back 
to von~Ettingshausen~\cite[Lec.\,42, Eq.\,(65)]{Ettingshausen:1827} in 1827.

Moreover, the results of \cite[Thms.\,1.2,\,1.3]{Kellner:2022} provide 
a denominator formula for $\BB_{r,\ell}$ implying that these numbers never vanish, 
except for the cases when $(r,\ell)$ equals $(n,0)$ or $(0,n)$ for odd $n \geq 3$.
An analogue for the Bernoulli polynomials is given by
\[
  \BB_{r,\ell}(x) = \sum_{\nu=0}^{r} \binom{r}{\nu} \BN_{\ell+\nu}(x),
\]
where $\BB_{r,\ell} = \BB_{r,\ell}(0)$. The reciprocity relation then holds 
asymmetrically by
\begin{equation} \label{eq:bpid}
  (-1)^r \, \BB_{r,\ell}(x) = (-1)^\ell \, \BB_{\ell,r}(-x).
\end{equation}
See Kellner~\cite[Thm.\,3.1]{Kellner:2021} for a short proof of \eqref{eq:bpid} 
using umbral calculus, but given more generally for Appell polynomials (see below), 
hence implying \eqref{eq:bpid} and so \eqref{eq:bnid} as special cases.

\begin{corollary} 
Let $n \geq k \geq 0$ and $\ell \geq 0$, then
\begin{align*}
  \sum_{\nu=0}^{k} \binom{n+1}{k-\nu} \BN_{\ell+\nu}(x)
  &= \sum_{\nu=0}^{k} \binom{n-\nu}{k-\nu} \BB_{\nu,\ell}(x)
\shortintertext{and}
  \sum_{\nu=0}^{k} \binom{n+1}{k-\nu} \BN_{\ell+\nu}
  &= \sum_{\nu=0}^{k} \binom{n-\nu}{k-\nu}\BB_{\nu,\ell}.
\shortintertext{In particular,}
  \sum_{\nu=0}^{n} \binom{n+1}{\nu+1} \BN_{\ell+\nu}(x)
  &= \sum_{\nu=0}^{n} \BB_{\nu,\ell}(x)
\shortintertext{and}
  \sum_{\nu=0}^{n} \binom{n+1}{\nu+1} \BN_{\ell+\nu}
  &= \sum_{\nu=0}^{n} \BB_{\nu,\ell}. 
\end{align*}
Moreover, $\BB_{\nu,\ell}(x)$ and $\BB_{\nu,\ell}$ can be replaced by
$(-1)^{\ell+\nu} \, \BB_{\ell,\nu}(-x)$ and $(-1)^{\ell+\nu} \, \BB_{\ell,\nu}$, 
respectively.
\end{corollary}

\subsection{Appell polynomials}

Let $(\AP_\nu(x))_{\nu \geq 0}$ be an Appell sequence of Appell polynomials 
$\AP_\nu(x)$ (see Appell~\cite{Appell:1880}). We set $\AP_\nu = \AP_\nu(0)$.
For $r,\ell \geq 0$, define the polynomials
\[
  \AB_{r,\ell}(x) = \sum_{\nu=0}^{r} \binom{r}{\nu} \AP_{\ell+\nu}(x)
\]
and the numbers $\AB_{r,\ell} = \AB_{r,\ell}(0)$.
If the Appell polynomials satisfy the reflection relation 
\begin{align}
  \AP_\nu(1-x) &= (-1)^\nu \AP_\nu(x),
  \label{eq:apref}
\shortintertext{then the reciprocity relation}
  (-1)^r \, \AB_{r,\ell}(x) &= (-1)^\ell \, \AB_{\ell,r}(-x)
  \nonumber
\shortintertext{holds (see \cite[Thm.\,3.1]{Kellner:2021}), as well as,}
  (-1)^r \, \AB_{r,\ell} &= (-1)^\ell \, \AB_{\ell,r}.
  \nonumber
\end{align}
For example, the Bernoulli and Euler polynomials $\BN_n(x)$ and $\EN_n(x)$ are 
Appell polynomials and satisfy the reflection relation \eqref{eq:apref}.

\begin{corollary} 
Let $n \geq k \geq 0$ and $\ell \geq 0$, then
\begin{align*}
  \sum_{\nu=0}^{k} \binom{n+1}{k-\nu} \AP_{\ell+\nu}(x)
  &= \sum_{\nu=0}^{k} \binom{n-\nu}{k-\nu} \AB_{\nu,\ell}(x)
\shortintertext{and}
  \sum_{\nu=0}^{k} \binom{n+1}{k-\nu} \AP_{\ell+\nu}
  &= \sum_{\nu=0}^{k} \binom{n-\nu}{k-\nu}\AB_{\nu,\ell}.
\shortintertext{In particular,}
  \sum_{\nu=0}^{n} \binom{n+1}{\nu+1} \AP_{\ell+\nu}(x)
  &= \sum_{\nu=0}^{n} \AB_{\nu,\ell}(x)
\shortintertext{and}
  \sum_{\nu=0}^{n} \binom{n+1}{\nu+1} \AP_{\ell+\nu}
  &= \sum_{\nu=0}^{n} \AB_{\nu,\ell}. 
\end{align*}
Moreover, if the reflection relation \eqref{eq:apref} holds, then
$\AB_{\nu,\ell}$ and $\AB_{\nu,\ell}(x)$ can be replaced by
$(-1)^{\ell+\nu} \AB_{\ell,\nu}$ and $(-1)^{\ell+\nu} \AB_{\ell,\nu}(-x)$, 
respectively.
\end{corollary}

\subsection{Fibonacci and Lucas numbers}

The sequences of the Fibonacci numbers (\seqnum{A000045}) and Lucas numbers 
(\seqnum{A000032}) for $n \geq 0$ begin
\begin{align*}
  \FN_n &= 0, 1, 1, 2, 3, 5, 8, 13, 21, 34, 55, 89, 144, 233, \ldots
\shortintertext{and}
  \LN_n &= 2, 1, 3, 4, 7, 11, 18, 29, 47, 76, 123, 199, 322, \ldots,
\end{align*}
obeying the same linear recurrences ${\FN_n = \FN_{n-1} + \FN_{n-2}}$ and 
${\LN_n = \LN_{n-1} + \LN_{n-2}}$ for ${n \geq 2}$, respectively. 
Since (cf.\ Hoggatt~\cite[Chap.\,9, p.\,61]{Hoggatt:1969})
\[
  \FN_{2n} = \sum_{\nu=0}^{n} \binom{n}{\nu} \FN_\nu \andq
  \LN_{2n} = \sum_{\nu=0}^{n} \binom{n}{\nu} \LN_\nu,
\]
these numbers satisfy the binomial transform. We obtain the following result 
due to linearity.

\begin{corollary} 
Let $\alpha, \beta \in \CC$. For $n \geq k \geq 0$, we have
\begin{align*}
  \sum_{\nu=0}^{k} \binom{n+1}{k-\nu} (\alpha \FN_\nu + \beta \LN_\nu)
  &= \sum_{\nu=0}^{k} \binom{n-\nu}{k-\nu} (\alpha \FN_{2\nu} + \beta \LN_{2\nu}).
\shortintertext{In particular,}
  \sum_{\nu=0}^{n} \binom{n+1}{\nu+1} (\alpha \FN_\nu + \beta \LN_\nu)
  &= \sum_{\nu=0}^{n} (\alpha \FN_{2\nu} + \beta \LN_{2\nu}).
\end{align*}
\end{corollary}

\subsection{Fibonacci matrix}

Let $R = \Mat_2(\ZZ)$. Define the Fibonacci matrix
\[
  \FM = \begin{pmatrix} 1 & 1 \\ 1 & 0 \end{pmatrix} \in R.
\]
Then we have the well-known identity (cf.~\cite[Chap.\,11, pp.\,65--67]{Hoggatt:1969})
\[
  \FM^\nu = \begin{pmatrix} \FN_{\nu+1} & \FN_{\nu} \\ \FN_{\nu} & \FN_{\nu-1} \end{pmatrix}
  \quad (\nu \geq 0),
\]
where $\FM^0 = \IM$ is the identity matrix and coincides with $\FN_{-1} = 1$ 
(extended by the recurrence relation ${\FN_1 = \FN_0 + \FN_{-1}}$).
Thus, it easily follows that $\FM + \IM = \FM^2$. Since $\IM\FM = \FM\IM = \FM$, 
this implies the binomial transform
\[
  \FM^{2n} = (\FM + \IM)^n = \sum_{\nu=0}^{n} \binom{n}{\nu} \FM^\nu.
\]
Finally, we obtain by Theorem~\ref{thm:main} for $R$ as above the following result.

\begin{corollary} 
For $n \geq k \geq 0$, we have
\begin{align*}
  \sum_{\nu=0}^{k} \binom{n+1}{k-\nu} \FM^\nu
  &= \sum_{\nu=0}^{k} \binom{n-\nu}{k-\nu} \FM^{2\nu}
\shortintertext{and}
  \sum_{\nu=0}^{k} \binom{n+1}{k-\nu} 
  \begin{pmatrix} \FN_{\nu+1} & \FN_{\nu} \\ \FN_{\nu} & \FN_{\nu-1} \end{pmatrix}
  &= \sum_{\nu=0}^{k} \binom{n-\nu}{k-\nu}
  \begin{pmatrix} \FN_{2\nu+1} & \FN_{2\nu} \\ \FN_{2\nu} & \FN_{2\nu-1} \end{pmatrix}.
\shortintertext{In particular,}
  \sum_{\nu=0}^{n} \binom{n+1}{\nu+1} \FM^\nu
  &= \sum_{\nu=0}^{n} \FM^{2\nu}
\shortintertext{and}
  \sum_{\nu=0}^{n} \binom{n+1}{\nu+1} 
  \begin{pmatrix} \FN_{\nu+1} & \FN_{\nu} \\ \FN_{\nu} & \FN_{\nu-1} \end{pmatrix}
  &= \sum_{\nu=0}^{n} 
  \begin{pmatrix} \FN_{2\nu+1} & \FN_{2\nu} \\ \FN_{2\nu} & \FN_{2\nu-1} \end{pmatrix}.
\end{align*}
\end{corollary}

\section{Proofs of main results}

As a key ingredient, we need the \emph{upside-down} convolution 
(cf.~\cite[Sec.\,5.1, Eq.\,(5.26), p.\,169]{GKP:1994}) 
for nonnegative integers $n$, $x$, and $y$ that
\begin{equation} \label{eq:conv}
  \sum_{\nu=0}^{n} \binom{\nu}{x} \binom{n-\nu}{y} = \binom{n+1}{x+y+1}.
\end{equation}

\begin{proof}[Proof of Theorem \ref{thm:main}]
Let $n \geq k \geq 0$.
Since both sides of \eqref{eq:seq4} are inverse to each other, we only
need to show the left-hand side of \eqref{eq:seq4}.
By symmetry of the binomial coefficients, \eqref{eq:seq1} and \eqref{eq:seq2}
are equivalent with a suitable change of indices. Now, \eqref{eq:seq3} is
a limiting case of \eqref{eq:seq1} by setting ${k=n}$, as well as, 
of \eqref{eq:seq2} by setting $k=0$. 

\phrase{$\Rightarrow$}:  
Note that ${s_0 = t_0}$. Applying the operator $\Delta$ to both sides 
of \eqref{eq:seq3} for $n \geq 1$, we infer that
\[
  \nabla^n \, s_\nu \valueat{\nu=0} =
  \Delta \sum_{\nu=0}^{n-1} \binom{n}{\nu+1} s_\nu \valueat{n}
  = \Delta \sum_{\nu=0}^{n-1} t_\nu \valueat{n} = t_n, 
\]
implying the result.

\phrase{$\Leftarrow$}:  
Expanding $t_\nu$ by the left-hand side of \eqref{eq:seq4}, 
the right-hand side of \eqref{eq:seq1} becomes
\begin{equation} \label{eq:R1}
  \rhs = \sum_{\nu=0}^{k} \binom{n-\nu}{n-k} \sum_{\ell=0}^{\nu} \binom{\nu}{\ell} s_\ell.
\end{equation}
This can be rewritten and rearranged to
\begin{equation} \label{eq:R2}
  \rhs = \sum_{\ell=0}^{k} \sum_{\nu=0}^{n} \binom{n-\nu}{n-k} \binom{\nu}{\ell} s_\ell
\end{equation}
by increasing the limits of the outer and inner sum of \eqref{eq:R1},
since the corresponding binomial coefficients $\binom{n-\nu}{n-k}$ and 
$\binom{\nu}{\ell}$ then vanish, respectively. 
The inner sum of \eqref{eq:R2} can be evaluated separately by
using the convolution identity \eqref{eq:conv}. Thus, we obtain
\[
  \sum_{\nu=0}^{n} \binom{n-\nu}{n-k} \binom{\nu}{\ell} 
  = \binom{n+1}{n-k+\ell+1} = \binom{n+1}{k-\ell}.
\]
This finally yields 
\[
  \rhs = \sum_{\ell=0}^{k} \binom{n+1}{k-\ell} s_\ell,
\]
being identical to the left-hand side of \eqref{eq:seq1}. 
This completes the proof.
\end{proof}

\begin{remark}
A similar proof was given in \cite{Kellner:2026} in connection with 
the Appell polynomials.  
\end{remark}

\begin{proof}[Proof of Corollary \ref{cor:hock}]
Let ${n \geq k \geq 0}$. We use Theorem~\ref{thm:main} and show two parts.

(1) Let $\seq{s} = (1,0,0,\ldots)$, where its binomial transform
is then given by $\seq{t} = (1,1,1,\ldots)$. From \eqref{eq:seq1}, we deduce
\[
  \binom{n+1}{k}
  = \sum_{\nu=0}^{k} \binom{n-\nu}{k-\nu} 
  = \sum_{\nu=0}^{k} \binom{n-k+\nu}{\nu},
\]
providing \eqref{eq:hs1}. Further, \eqref{eq:seq2} directly implies \eqref{eq:hs2}.

(2) Let $\seq{s} = (1,-1,1,-1,\ldots)$. By $(1-1)^\ell = 0^\ell$ for 
$\ell \geq 0$, the binomial transform is given by $\seq{t} = (1,0,0,\ldots)$. 
By \eqref{eq:seq1}, we have
\[
  \sum_{\nu=0}^{k} \binom{n+1}{k-\nu} (-1)^\nu
  = \sum_{\nu=0}^{k} \binom{n+1}{\nu} (-1)^{k-\nu}
  = \binom{n}{k},
\]
which is \eqref{eq:hs3}. Finally, \eqref{eq:seq2} instantly shows \eqref{eq:hs4},
completing the proof.
\end{proof}

\begin{proof}[Proof of Theorem \ref{thm:func}]
Let $f: \RR \to \CC$. We set $\seq{t} = \seqs{f(\nu)}$ and 
$\seq{s} = \seqs{\Delta^{\!\nu} t} = \seqs{\Delta^{\!\nu} f}$.
By Theorem~\ref{thm:main}, the relations \eqref{eq:fp1} -- \eqref{eq:fp3} then 
are equivalent to \eqref{eq:seq1} -- \eqref{eq:seq3}, respectively.
Now, let $f \in \CC[x]$ be a polynomial with $d = \deg f$. If $f = 0$, then all
equations are empty, and we are done. Otherwise, we have that $d \geq 0$ and
the sum of the left-hand side of \eqref{eq:fp1} -- \eqref{eq:fp3}, respectively, 
becomes truncated when $\nu > d$, since then $\Delta^{\!\nu} f = 0$.
Clearly, if $f \in \Int(\ZZ)$, then $\seqs{\Delta^{\!\nu} f}$ lies in $\ZZ$, and
so \eqref{eq:fp1} -- \eqref{eq:fp3} also lie in $\ZZ$. This completes the proof.
\end{proof}

\begin{proof}[Proof of Corollary \ref{cor:bin}]
Let $d \geq 0$ and define ${f(\nu) = \binom{x+\nu}{d}}$ for $\nu \geq 0$.
Note that $\Delta \binom{x}{k} = \binom{x}{k-1}$. As a consequence,
\[
  \Delta^{\!\nu} f = \Delta^{\!\nu} \binom{x}{d} = \binom{x}{d-\nu}
  \quad (\nu \geq 0).
\]
Applying Theorem \ref{thm:func} shows the result, completing the proof.
\end{proof}

\appendix
\section{Summary}

Let $n \geq k \geq 0$. We summarize the obtained boomerang identities of the paper by
\begin{align*}
  \sum_{\nu=0}^{k} \binom{n+1}{k-\nu} \, \RT_\nu
  &= \sum_{\nu=0}^{k} \binom{n-\nu}{k-\nu} \, \CT_\nu \\
\shortintertext{and}  
  \sum_{\nu=k}^{n} \binom{n+1}{\nu+1} \RT_{\nu-k} 
  &= \sum_{\nu=k}^{n} \binom{\nu}{k} \CT_{n-\nu},
\shortintertext{as well as the limiting case}
  \sum_{\nu=0}^{n} \binom{n+1}{\nu+1} \, \RT_\nu
  &= \sum_{\nu=0}^{n} \CT_\nu, 
\end{align*}
where $\RT$ and $\CT$ represent terms in the rows and columns of Pascal's triangle, 
respectively. In the following table, let $d, \ell, m \geq 0$ and $x$ as an unknown. 
Let $\nu \geq 0$, except for the cases $\frac{1}{\nu}$ and $H_\nu$ where $\nu \neq 0$.
\medskip

\begingroup
\small
\setstretch{1.3}
\begin{longtable}{ccc}
  \toprule
  $\RT_\nu$ & $\CT_\nu$ & Application \\
  \midrule
  \endhead
  
  \multicolumn{3}{r}{\footnotesize\textit{Continued on next page}} \\
  \endfoot

  \bottomrule
  \addlinespace[6pt]
  \caption{Binomial transforms.} \label{tab:sum} \\
  \endlastfoot
  
  $\seq{s}$ & $\seqs{\nabla^\nu\! s}$ & Sequence \\
  $\seqs{\Delta^{\!\nu} t}$ & $\seq{t}$ & \ditto \\
  \midrule
  $\Delta^{\!\nu} f$ & $f(\nu)$ & Function \\
  $\binom{x}{d-\nu}$ & $\binom{x+\nu}{d}$ & Binomial \\
  \midrule
  $(1,0,0,\ldots)$ & $(1,1,1,\ldots)$ & Hockey stick \\
  $(1,-1,1,-1,\ldots)$ & $(1,0,0,\ldots)$ & \ditto \\
  \midrule
  $\BE_\nu$ & $\BE_{\nu+1}$ & Bell number \\
  $x \, \BE_\nu(x)$ & $\BE_{\nu+1}(x)$ & Bell polynomial \\
  \midrule
  $\nu! \, \ST_2(m,\nu)$ & $\nu^m$ & Stirling number $\ST_2$ \\
  $(-1)^{\nu+1} \frac{1}{\nu}$ & $\HN_\nu$ & Harmonic number \\
  \midrule
  $\AP_{\ell+\nu}(x)$ & $\AB_{\nu,\ell}(x)$ & Appell polynomial \\
  $\BN_{\ell+\nu}$ & $\BB_{\nu,\ell}$ & Bernoulli number \\
  $\BN_{\ell+\nu}(x)$ & $\BB_{\nu,\ell}(x)$ & Bernoulli polynomial \\
  \midrule
  $\FN_\nu$ & $\FN_{2\nu}$ & Fibonacci number \\
  $\LN_\nu$ & $\LN_{2\nu}$ & Lucas number \\
  \midrule
  $\begin{pmatrix} 1 & 1 \\ 1 & 0 \end{pmatrix}^{\!\!\nu}$ & 
  $\begin{pmatrix} 1 & 1 \\ 1 & 0 \end{pmatrix}^{\!\!2\nu}$ &
  Fibonacci matrix \\
  $\begin{pmatrix} \FN_{\nu+1} & \FN_{\nu} \\ \FN_{\nu} & \FN_{\nu-1} \end{pmatrix}$ &
  $\begin{pmatrix} \FN_{2\nu+1} & \FN_{2\nu} \\ \FN_{2\nu} & \FN_{2\nu-1} \end{pmatrix}$ &
  \ditto \\
\end{longtable}
\endgroup

\bibliographystyle{amsplain}

\end{document}